\documentclass[12pt]{amsart}
\usepackage{amssymb,amsmath}
\DeclareMathOperator{\dvg}{div}
\newcommand{\N}{{\mathbb N}}

\newcommand{\R}{{\mathbb R}}
\newcommand{\eps}{\varepsilon}
\newcommand{\de}{\partial}
\renewcommand{\theta}{\vartheta}

\numberwithin{equation}{section}
\newtheorem{theorem}{Theorem}[section]
\newtheorem{proposition}[theorem]{Proposition}
\newtheorem{corollary}[theorem]{Corollary}
\newtheorem{lemma}[theorem]{Lemma}
\newtheorem{definition}[theorem]{Definition}
\theoremstyle{definition}
\newtheorem{remark}[theorem]{Remark}

\title[Location of spikes
for quasi-linear elliptic equations]{On
the location of concentration points for
singularly perturbed elliptic equations}

\author{Simone Secchi}
\author{Marco Squassina}

\address{Dipartimento di Matematica ``L.Tonelli''
\newline\indent
Universit\`a degli Studi di Pisa
\newline\indent
Via Buonarroti 2, I-56127 Pisa, Italy}
\email{secchi@dm.unipi.it}

\address{Dipartimento di Matematica ``F.Brioschi''
\newline\indent
Politecnico di Milano
\newline\indent
Via Bonardi 9, I-20133 Milano, Italy}
\email{squassina@mate.polimi.it}

\thanks{The first author was supported by the MIUR national research
project ``Variational Methods and Nonlinear Differential Equations''.
The second author was supported by the MIUR national research
project ``Variational and Topological
Methods in the Study of Nonlinear Phenomena'' and
by the Istituto Nazionale di Alta Matematica ``F. Severi'' (INdAM)}


\begin{document}

\begin{abstract}
By exploiting a variational identity of Poho\v zaev-Pucci-Serrin type
for solutions of class $C^1$, we get
some necessary conditions for locating the peak-points of a
class of singularly perturbed quasi-linear elliptic problems in divergence
form. More precisely, we show that the points where
the concentration occurs, in general, must belong to what we call
the set of weak-concentration points.
Finally, in the semi-linear case, we provide a new necessary
condition which involves the Clarke subdifferential
of the ground-state function.
\end{abstract}
\maketitle


\section{Introduction}

\noindent
Let $\eps >0$, $n \geq 3$, and $1<p<n$.
In this paper we consider
the following class of singularly perturbed quasi-linear elliptic
problems in divergence form
\begin{equation}
\label{problem}
\tag{$P_\eps$}
\begin{cases}
-\eps^p\dvg(\alpha(x)\nabla\beta(\nabla u))
+V(x)u^{p-1}=K(x)f(u) & \text{in $\R^n$}\\
u>0 & \text{in $\R^n$}.
\end{cases}
\end{equation}
We assume that the functions $\alpha$, $V$, $K \colon \R^n\to\R$
are positive, of class $C^1$ with bounded
derivatives, and $\alpha$, $K\in L^\infty (\R^n)$. Moreover, let
\begin{equation*}
\inf_{x\in\R^n}\alpha(x)>0
\quad\text{and}\quad
\inf_{x\in\R^n}V(x)>0.
\end{equation*}
The function $\beta \colon \R^n\to\R$ is of class $C^1$, strictly convex, and
positively homogeneous of degree $p$, namely,
$\beta (\lambda \xi) = \lambda^p \beta (\xi)$
for every $\lambda >0$ and $\xi\in\R^n$. Moreover, there exist $\nu>0$ and
$c_1,c_2>0$ such that
\begin{gather}
\label{betacont}
\nu|\xi|^p\leq \beta(\xi)\leq c_1|\xi|^p, \\
\label{betacontxi}
|\nabla \beta(\xi)|\leq c_2|\xi|^{p-1},
\end{gather}
for every $\xi\in\R^n$. The nonlinearity $f:\R^+\to\R$
is of class $C^1$ and such that
\begin{equation*}
\lim_{s\to 0^+}\frac{f(s)}{s^{p-1}}=0
\quad\text{and}\quad
\lim_{s\to+\infty}\frac{f(s)}{s^{q-1}}=0,
\end{equation*}
for some $p<q<p^*$, with $p^*=np/(n-p)$. Moreover,
\begin{equation*}
0<\vartheta F(s)\leq f(s)s,\quad\text{for every $s>0$},
\end{equation*}
for some $\vartheta>p$, where we have set
$F(s)=\int_0^s f(t)\,dt$, $s\in\R^+$.
\par
Let us define the space $W_V(\R^n)$ by setting
\[
W_V(\R^n) :=\Big\{ u\in W^{1,p}(\R^n):\,\,
\int_{\R^n} V(x) |u|^p\,dx < \infty \Big\},
\]
endowed with the natural norm
$\|u\|_{W_V}^p=\int_{\R^n}|\nabla u|^p\,dx+\int_{\R^n} V(x)|u|^p\,dx$.
For $p=2$, we write $H_V(\R^n)$ in place of $W_V(\R^n)$.
Under the previous assumptions, if $K\equiv 1$, it has been
recently proved in~\cite{giacsqu} (see also~\cite{squa}) that if for some compact
subset $\Lambda\subset\R^n$ we have
\begin{equation*}
V(z_0)=\min_{\Lambda} V <\min_{z\in\partial\Lambda} V(z)
\quad\text{and}\quad
\alpha(z_0)=\min_{z\in\Lambda}\alpha(z),
\end{equation*}
then, for every $\eps$ sufficiently small, there exists
a solution $u_\eps\in W_V(\R^n)$ of~\eqref{problem}
which has a maximum point $z_{\eps}\in\Lambda$, with
\begin{equation*}
\lim_{\eps\to 0}V(z_\eps)=\min_{\Lambda} V
\quad\text{and}\quad
\lim_{\eps\to 0}\|u_\eps\|_{L^{\infty}(\Omega\setminus B_\rho(z_{\eps}))}=0,
\quad\text{for every $\rho>0$}.
\end{equation*}
In the semi-linear case, the construction of solutions concentrating
at critical points (or minima) of the potential $V(x)$ or other finite
dimensional driven functions has been deeply investigated in the last decade
and also stronger results can be found in the literature
(see e.g.\ \cite{abc,dancer,pinofelm,pinofelm3,pinofelm4,
pinofelm5,oh1,rabinowitz,wang} and references therein).
\par
The goal of this paper is to establish some {\em necessary conditions} for a
sequence of solutions $(u_{\eps_h})$ of~\eqref{problem} to concentrate
around a given point $z_0\in\R^n$, in the sense of
Definition~\ref{concentrationdef}.
If $\beta(\xi)=\xi$, we will prove (see Theorem~\ref{laplaci})
that if $z_0$ is a concentration point for a sequence $(u_{\eps_h})\subset H_V(\R^n)$
of solutions of the problem then there exists a locally Lipschitz function
$\Sigma:\R^n\to\R$, the ground-state
function, which has, under suitable assumptions, a critical point
in the sense of the Clarke subdifferential at $z_0$, that is,
\begin{equation*}
0\in\partial\Sigma(z_0).
\end{equation*}
Under more stringent assumptions, it turns out that $\Sigma$
admits all the directional derivatives at $z_0$ and
$\nabla \Sigma(z_0)=0$.
In the general case, as a first necessary
condition, the gradient vectors
\begin{equation*}
\nabla\alpha(z_0),\,\,\,
\nabla V(z_0),\,\,\,
\nabla K(z_0)
\end{equation*}
must be linearly dependent. Moreover,
in Theorem~\ref{mainth} (see also Theorem~\ref{plaplaci}),
we show that the concentration points for problem~\eqref{problem}
must belong to a set $\mathfrak{C}$ (which has a variational structure) that
we call the set of weak-concentration points (see Definition~\ref{gammapm}).
To the authors' knowledge, this kind of necessary
conditions in terms of generalized
gradients seem to be new. Quite interestingly,
the lack of uniqueness (up to translations) for the limiting problem
(namely the rescaled problem with frozen coefficients)
\begin{equation}
\label{limprob}
\tag{$P_{z}$}
\begin{cases}
-\alpha(z)\dvg(\nabla\beta(\nabla u))
+V(z)u^{p-1}=K(z)f(u) & \text{in $\R^n$}\\
u>0 & \text{in $\R^n$}
\end{cases}
\end{equation}
induces a lack of regularity for $\Sigma$.
Some conditions ensuring uniqueness of solutions for~\eqref{limprob}
can be found in~\cite{dr,sertang}.
For instance, for $1<p\leq 2$, $\beta (\xi)=|\xi|^{p-2}\xi$,
and $f(u)= u^{q-1}$ with $p<q<p^*$, we have uniqueness
and $\Sigma$ admits all the directional derivatives.
\par
We stress that some necessary conditions for the
location of concentration points were previously obtained
by Ambrosetti et al.\ in~\cite{abc}, and by Wang and Zeng
in~\cite{wang,wangzeng} in the case $p=2$
and $\beta(\xi)=\xi$. Their approach is
based on a repeated use of the Divergence Theorem.
With respect to those
papers we prove our main results
by means of a locally Lipschitz variant of the celebrated Pucci-Serrin
variational identity \cite{ps}. In our possibly degenerate
setting, classical $C^2$ solutions might not exist, the highest general
regularity class being $C^{1,\beta}$ (see \cite{tolks}). Therefore, the
classical identity is not applicable in our framework. However, it
has been recently shown in \cite{dms} that,
under minimal regularity assumptions, the identity
holds for locally Lipschitz solutions (see Theorem~\ref{step1}), provided that
the operator ($\beta$, in our case) is strictly convex in the
gradient, which, from our viewpoint, is a
very natural requirement.
\par
This identity has also turned out to be useful in
characterizing the exact energy level of the least
energy solutions of the problem~\eqref{limprob}.
Indeed in~\cite[Theorem 3.2]{giacsqu} it was proved that~\eqref{limprob}
admits a least energy solution $u_z\in W^{1,p}(\R^n)$ having the
Mountain-Pass energy level. This is precisely the motivation
that led us to define the ground-state function
$\Sigma$ also in a degenerate setting.

\section{The Quasi-linear Case}

The aim of this section is the study of some necessary
conditions for the concentration of the solutions at a point $z_0$
to occur, in the quasi-linear framework.

\subsection{Some preliminary definitions and properties}
If $z$ is fixed in $\R^n$, we consider the limiting
functional $I_z:W^{1,p}(\R^n)\to\R$,
\begin{equation*}
I_{z}(u):=\alpha(z)\int_{\R^n} \beta(\nabla u)\,dx
+\frac{V(z)}{p}\int_{\R^n}|u|^p\,dx-K(z)\int_{\R^n}F(u)\,dx.
\end{equation*}
It follows from our assumptions on $\beta$ and $f$ that $I_z$ is a $C^1$
functional, and its critical points are solutions of the
limiting problem~\eqref{limprob}. We
define the minimax value $c_z$ for $I_z$ by setting
\begin{align}
\label{MPlev}
c_{z}&:=\inf_{\gamma\in {\mathcal P}_{z}}
\sup_{t\in[0,1]}I_{z}(\gamma(t)),\\
{\mathcal P}_{z} &:=\Big\{\gamma\in
C([0,1],W^{1,p}(\R^n)):\,\,\gamma(0)=0,\,\,\,
I_{z}(\gamma(1))<0\Big\}.\notag
\end{align}
Throughout the rest of the paper, we will denote by $G(z)$ the set of
all the nontrivial solutions, up to translations, of the limiting problem~$(P_z)$
(the set of bound-states). Under our assumptions
on $f$, $G(z)\not=\emptyset$ for every $z\in\R^n$. Finally, $\cdot$
will always stand for the usual inner product of $\R^n$.
\vskip2pt
We now introduce two functions $\partial\Gamma^-,\partial\Gamma^+$
that will be useful in the sequel.

\begin{definition}
\label{gammapm}
For every $z,w\in\R^n$ we define $\partial\Gamma^-(z;w)$
and $\partial\Gamma^+(z;w)$ by setting
\begin{align*}
&\partial\Gamma^-(z;w):=\sup_{v\in G(z)}\nabla_z I_z(v)\cdot w, \\
&\partial\Gamma^+(z;w):=\inf_{v\in G(z)}\nabla_z I_z(v)\cdot w,
\end{align*}
where $\nabla_z$ denotes the gradient with respect to $z$.
Explicitly, for every $z,w\in\R^n$,
\begin{align*}
\partial\Gamma^-(z;w)&=
\sup_{v\in G(z)}\bigg[\nabla \alpha(z)\cdot w\int_{\R^n}
\beta(\nabla v)\,dx\\
&+\nabla V(z)\cdot w\int_{\R^n}\frac{|v|^p}{p}\,dx-
\nabla K(z)\cdot w\int_{\R^n}F(v)\,dx\bigg],\\
\partial\Gamma^+(z;w)&=
\inf_{v\in G(z)}\bigg[\nabla \alpha(z)\cdot w\int_{\R^n}
\beta(\nabla v)\,dx\\
&+\nabla V(z)\cdot w\int_{\R^n}\frac{|v|^p}{p}\,dx
-\nabla K(z)\cdot w\int_{\R^n}F(v)\,dx\bigg],
\end{align*}
Finally, we define a set $\mathfrak{C}\subset\R^n$ by
\begin{equation*}
\mathfrak{C}:=
\Big\{z\in\R^n:\,\,\,\text{$\partial\Gamma^-(z,w)\geq 0$
and\, $\partial\Gamma^+(z,w)\leq 0$,
for every $w\in\R^n$}\Big\}.
\end{equation*}
We say that $\mathfrak{C}$ is the set of weak-concentration points
for problem~\eqref{problem}.
\end{definition}
\noindent
The motivations that lead us to introduce the functions $\partial\Gamma^-$,
$\partial\Gamma^+$ and the set of weak-concentration points
will be clear in the course of the investigation.
\vskip2pt
For the sake of completeness, we recall the following

\begin{definition}
\label{groundef}
We define the {\em ground-state function}
$\Sigma:\R^n\to\R$ by setting
\begin{equation*}
\Sigma(z):=\min_{u\in G(z)}I_z(u),\quad\text{for every $z\in\R^n$.}
\end{equation*}
\end{definition}

We now collect a few useful properties of the function $\Sigma$.
\begin{lemma}
\label{pro-sigma}
Assume that
\begin{equation}
\label{increasing}
\text{the map\,\, $s\in \R^+ \mapsto
\frac{f(s)}{s^{p-1}}$\,\, is increasing.}
\end{equation}
Then, the following facts hold:
\begin{itemize}
\item[(i)] the map $\Sigma$ is well defined,
continuous, and
\begin{equation*}
\Sigma(z)=c_z,\quad \text{for every $z\in\R^n$};
\end{equation*}
\vskip1pt
\item[(ii)] the map $\Sigma$ can be written as
\begin{equation*}
\Sigma (z)=\!\!\!\!\inf_{u\in W^{1,p}(\R^n)\setminus\{0\}}
\max_{\theta\geq 0}I_z(\theta u)
=\inf_{u\in \mathcal{N}_z} I_z (u),\quad
\text{for every $z\in\R^n$},
\end{equation*}
where $\mathcal{N}_z$ is the Nehari manifold, defined as
\begin{equation*}
\mathcal{N}_z:=\Big\{ u\in W^{1,p}(\R^n) \setminus \{0\}:\,\,
I_z'(u)[u]=0 \Big\}.
\end{equation*}
\end{itemize}
\end{lemma}
\begin{proof}
To prove {\rm (ii)}, it suffices to argue as in
\cite[Proposition 2.5]{pompsec}. We now come to assertion {\rm (i)}.
By~\cite[Theorem 3.2]{giacsqu}, for every $z\in \R^n$,
problem~\eqref{limprob} admits a solution
$v_z\in W^{1,p}(\R^n)$, $v_z\neq 0$, such that
$$
I_{z}(v_z)=\Sigma(z)=c_{z},
$$
where $c_z$ is defined as in~\eqref{MPlev}.
The continuity of $\Sigma$ then follows from the continuity
of the map $z\mapsto c_z$, which we now prove directly
using an argument envisaged by Rabinowitz~\cite{rabinowitz}.
For $\alpha$, $V$, $K\in \R$, define the
functional $I_{\alpha,V,K}:W^{1,p}(\R^n) \to\R$ by
\begin{equation*}
I_{\alpha,V,K} (u):= \alpha \int_{\R^n}
\beta(\nabla u)\, dx + \frac{V}{p}\int_{\R^n}
|u|^p \, dx - K \int_{\R^n} F(u) \, dx.
\end{equation*}
Let us set:
\begin{align*}
c(\alpha,V,K)&:=\inf_{\gamma\in\mathcal{P}_{
\alpha,V,K}} \max_{t\in [0,1]}
I_{\alpha,V,K} (\gamma (t)),\\
\noalign{\vskip2pt}
\mathcal{P}_{\alpha,V,K}&:=\Big\{ \gamma
\in C([0,1],W^{1,p}(\R^n)):\,\,
\gamma (0)=0,\,\, I_{\alpha,V,K}(\gamma(1)) < 0 \Big\}.
\end{align*}
\noindent
\emph{Claim}:
For every $(\alpha,V,K)\in\R^3$ we have
\begin{equation*}
\lim_{\eta \to 0} c(\alpha+\eta,V+\eta,K-\eta) = c(\alpha,V,K).
\end{equation*}
We first observe that a simple adaptation
of the argument of~\cite[Lemma 3.17]{rabinowitz} yields
\begin{equation}
\label{dismp}
\alpha_1>\alpha_2,\,\,
V_1>V_2,\,\, K_1<K_2\,\,\,\Longrightarrow\,\,\,
c(\alpha_1,V_1,K_1) \geq c(\alpha_2,V_2,K_2).
\end{equation}
The proof of the claim will be accomplished
indirectly. By virtue of~\eqref{dismp}, we get
\begin{equation*}
\lim_{\eta \to 0^-} c(\alpha+\eta,V+\eta,K-\eta):= c^- \leq c(\alpha,V,K).
\end{equation*}
Suppose that $c^- < c(\alpha,V,K)$. For the sake of brevity, we define
\begin{equation*}
J_\eta (u):=I_{\alpha+\eta,V+\eta,K-\eta}(u).
\end{equation*}
Let $\eta_h \to 0^-$ as $h\to\infty$, and $\delta_j \to 0^+$ as
$j\to\infty$. For each $h\in\N$, by assertion {\rm (ii)}, there is a sequence
$(u_{hj})$ in $W^{1,p}(\R^n)$, $u_{hj}\neq 0$, such that
\begin{equation}
\label{eq:equiv}
\alpha \int_{\R^n} \beta(\nabla u_{hj})\,dx+
V\int_{\R^n} |u_{hj}|^p\,dx =1
\end{equation}
and
\begin{equation}
\label{maxinequal}
\max_{\theta \geq 0} J_{\eta_h} (\theta u_{hj})
\leq c(\alpha+\eta_h,V+\eta_h,K-\eta_h) + \delta_j.
\end{equation}
Notice that we can choose the sequence
$(u_{hj})$ satisfying~\eqref{eq:equiv}, since the position
$$
u\mapsto \alpha \int_{\R^n} \beta(\nabla u)\,dx+V \int_{\R^n} |u|^p\,dx
$$
defines on $W^{1,p}(\R^n)$ a norm equivalent to
the natural one, as follows from~\eqref{betacont}.
Take now $j=h$, and set $u_h=u_{hh}$. Hence, in view
of~\eqref{maxinequal}, we have
\begin{align*}
c(\alpha,V,K) & \leq \max_{\theta \geq 0}
I_{\alpha,V,K}(\theta u_h) = I_{\alpha,V,K}(\phi(u_h)u_h) \\
&= J_{\eta_h} (\phi(u_h)u_h) - \eta_h \phi(u_h)^p
\int_{\R^n} \frac{|u_h|^p}{p}\,dx - \eta_h \phi(u_h)^p \int_{\R^n} \beta(\nabla u_h)\,dx \\
&\null \qquad - \eta_h \int_{\R^n} F(\phi(u_h)u_h)\,dx \\
&\leq \max_{\theta\geq 0} J_{\eta_h} (\theta u_h) -
\eta_h \phi(u_h)^p \int_{\R^n} \frac{|u_h|^p}{p}\,dx - \eta_h \phi(u_h)^p \int_{\R^n}
\beta(\nabla u_h)\,dx \\
&\null \qquad - \eta_h \int_{\R^n} F(\phi(u_h)u_h)\,dx \\
&\leq c(\alpha+\eta_h,V+\eta_h,K-\eta_h) + \delta_h -
\eta_h \phi(u_h)^p \int_{\R^n} \frac{|u_h|^p}{p}\,dx \\
&\null \quad - \eta_h \phi(u_h)^p \int_{\R^n} \beta(\nabla u_h)\,dx
+ \eta_h \int_{\R^n} F(\phi(u_h)u_h)\,dx \\
&\leq c^- + \delta_h - \eta_h \phi(u_h)^p
\int_{\R^n} \frac{|u_h|^p}{p}\,dx - \eta_h \phi(u_h)^p \int_{\R^n} \beta(\nabla u_h)\,dx \\
&\null \qquad - \eta_h \int_{\R^n} F(\phi(u_h)u_h)\,dx.
\end{align*}
At this point, one can show exactly as in
\cite[pp.~281-282]{rabinowitz} that there exists a
constant $C>0$ such that $\phi(u_h)\leq C$,
for every $h\in\N$ sufficiently large.
Therefore, recalling the properties of $F$ and the
Sobolev embedding, the above chain of inequalities contradicts
$c^- < c(\alpha,V,K)$, at least for every $h\in\N$
large enough. We conclude that $c^- < c(\alpha,V,K)$ is impossible.
In a completely similar fashion one can prove that the inequality
\begin{equation*}
c(\alpha,V,K) <\lim_{\eta\to 0^+} c(\alpha+\eta,V+\eta,K-\eta)
\end{equation*}
leads to a contradiction. Therefore the claim is proved.
\vskip1pt
\noindent
Let now $(z_h)$ be a sequence in $\R^n$ such that $z_h\to z$ as
$h\to\infty$. Observe that, given $\eta >0$,
for large $h\in\N$, we have
\begin{align*}
V(z)+\eta &\geq V(z)+|V(z_h)-V(z)| \\
&\geq V(z) \geq V(z) - |V(z_h)-V(z)| \\
&\geq V(z) - \eta,
\end{align*}
and similar relations hold for $\alpha$ and $K$.
Therefore the continuity of $z \mapsto c_z$
follows from the previous claim, applied
with $\alpha=\alpha (z)$, $V=V(z)$, and $K=K(z)$.
This completes the proof of assertion {\rm (i)}.
\end{proof}

\begin{remark}
As we have already pointed out in the introduction, we believe that
the lack of regularity of the ground-state map $\Sigma$ is essentially inherited by
the lack of uniqueness assumptions on the limiting equation~\eqref{limprob}.
From this viewpoint, in the degenerate case $p\neq 2$, the problem
of establishing the regularity of $\Sigma$ seems quite a difficult matter.
On the contrary, if $p=2$ and, for instance, $\beta(\xi)=\xi$, it is known
that $\Sigma$ is always at least locally Lipschitz continuous
(cf.\ Lemma~\ref{pro-sigma0}). If, additionally, $f(u)$ is exactly the power $u^{p-1}$
(in which case equation~\eqref{limprob} has in fact a unique solution~\cite{chen-lin}),
then $\Sigma$ is smooth and it also admits an explicit representation formula
(see Remark~\ref{smooth-power}).
\end{remark}

\noindent
Let now ${\mathcal L}:\R^n\times\R\times\R^n\to\R$
be a function of class $C^1$ such that
\begin{equation*}
\text{the function\,\,
$\xi\mapsto {\mathcal L}(x,s,\xi)$\,\,
is strictly convex},
\end{equation*}
for every $(x,s) \in \R^n \times \R$,
and let $\varphi \in L^{\infty}_{{\rm loc}}(\R^n)$.
\par
Next, we recall a Pucci-Serrin variational
identity for locally Lipschitz continuous solutions
of a general class of Euler equations, recently
proved in~\cite{dms}. As we have already
remarked in the introduction, the classical identity \cite{ps} is not applicable
here, since it requires the $C^2$ regularity of the solutions, while
the maximal regularity for degenerate equations
is $C^{1,\beta}$ (see e.g.\ \cite{tolks}).

\begin{theorem}
\label{step1}
Let $u:\R^n\to\R$ be a locally Lipschitz solution of
\begin{equation*}
- \dvg \!\left(\de_{\xi}{\mathcal L}(x,u,\nabla u)\right) +
\de_s{\mathcal L}(x,u,\nabla u)=\varphi\quad
\text{in ${\mathcal D}'(\R^n)$}.
\end{equation*}
\indent
Then,
\begin{multline}
\label{prima}
 \sum_{i,j=1}^n\int_{\R^n} \de_i h^j
\de_{\xi_i}{\mathcal L}(x,u,\nabla u) \de_j u\,dx  + \mbox{} \\
-
\int_{\R^n}\bigl[
(\dvg h) \, {\mathcal L}(x,u,\nabla u) +
h\cdot \de_{x}{\mathcal L}(x,u,\nabla u)
\bigr]\,dx =
\int_{\R^n} (h\cdot \nabla u)\varphi\,dx,
\end{multline}
for every $h\in C^1_c(\R^n,\R^n)$.
\end{theorem}

\subsection{Necessary conditions for locating peak-points}
We now state and prove the main results of this section.

\begin{theorem}
\label{mainth}
Let $z_0\in\R^n$, and assume that $(u_{\eps_h})$ is a
sequence of solutions of problem~\eqref{problem} such that
\begin{equation}
\label{concentration}
u_{\eps_h} = v_0\left(\frac{\cdot -z_0}{\eps_h}\right)+o(1),
\quad\text{strongly in $W_V(\R^n)$},
\end{equation}
for some $v_0\in W_V(\R^n)\setminus \{0\}$. Then, the following facts hold:
\begin{itemize}
\item[(a)] the vectors
\begin{equation*}
\nabla\alpha(z_0),\,\,\,
\nabla V(z_0),\,\,\,
\nabla K(z_0)
\end{equation*}
are linearly dependent\,;
\vskip3pt
\item[(b)] $z_0\in \mathfrak{C}$, that is
$z_0$ is a weak-concentration point for~\eqref{problem};
\vskip3pt
\item[(c)] if $G(z_0)=\{v_0\}$, then
all the partial derivatives of $\Sigma$ at $z_0$ exist, and
\begin{equation*}
\nabla\Sigma(z_0)=0,
\end{equation*}
that is $z_0$ is a critical point of $\Sigma$.
\end{itemize}
\end{theorem}
\begin{proof}
We write $u_h$ in place of $u_{\eps_h}$,
and we define
\begin{equation} \label{vh}
v_h(x):=u_h(z_0+\eps_hx).
\end{equation}
Therefore, $v_h$ satisfies the rescaled equation
\[
-\dvg(\alpha(z_0+\eps x) \nabla\beta
(\nabla v_h))+V(z_0+\eps x) v_h^{p-1}
=K(z_0+\eps x) f(v_h) \quad\text{in $\R^n$.}
\]
By \eqref{concentration}, we have $v_h \to
v_0$ strongly in $W_V(\R^n)$.
We now prove that $v_h\to v_0$ in the
$C^1$ sense over the compact sets of $\R^n$
and that $v_0$ is a nontrivial positive solution of
the equation
\begin{equation}
\label{sollim}
-\alpha(z_0){\rm div}(\nabla\beta(\nabla v))+
V(z_0)v^{p-1}=K(z_0)f(v)\quad\text{in $\R^n$}.
\end{equation}
Let us set
\begin{align*}
d_{h}(x)&:=
\begin{cases}
V(z_0+{\eps_h}x)-K(z_0+\eps_hx)\frac{f(v_h(x))}{v_h^{p-1}(x)} & \text{if
$v_h(x)\neq 0$} \\0 & \text{if $v_{h}(x)=0$},
\end{cases}\\
A(x,s,\xi)&:=\alpha(z_0+\eps_hx)\nabla \beta(\xi), \\
B(x,s,\xi)&:=d_{h}(x)s^{p-1},
\end{align*}
for every $x\in\R^n$, $s\in\R^+$, and $\xi\in\R^n$.
Taking into account~\eqref{betacontxi}, and the strict
convexity of $\beta$, we get
\begin{equation*}
A(x,s,\xi)\cdot \xi \geq\nu |\xi|^p
\quad\text{and}\quad
|A(x,s,\xi)|\leq c_2|\xi|^{p-1}.
\end{equation*}
Notice that, in view of the growth assumptions on $f$, there exists
$\delta > 0$ sufficiently small such that
$d_{h}\in L^{n/(p-\delta)}(B_{2\rho})$
for every $\rho>0$ and
\begin{equation*}
S=\sup_{h\in\N}\|d_{h}\|_{L^{n/(p-\delta)}
(B_{2\rho})}\leq D_\rho\Big(1+\sup_{h\in\N}
\|v_{h}\|_{L^{p^*}(B_{2\rho})}\Big)<\infty,
\end{equation*}
for some $D_\rho>0$.
Since we have ${\rm div}(A(x,v_h,\nabla v_h))=
B(x,v_h,\nabla v_h)$ for every $h\in\N$,
by exploiting~\cite[Theorem 1]{serrin}
there exists a radius $\rho>0$ and a
positive constant $M=M(\nu,c_2,S \rho^\delta)$ such that
\begin{equation*}
\sup_{h\in\N}\max_{x\in B_\rho}|v_{h}(x)|
\leq M(2\rho)^{-N/p}
\sup_{h\in\N}\|v_{h}\|_{L^p(B_{2\rho})}<\infty,
\end{equation*}
so that $(v_{h})$ is uniformly bounded in $B_\rho$.
Then, by virtue of~\cite[Theorem 8]{serrin},
up to a subsequence $(v_{h})$ converges uniformly
to $v_0$ in a small neighborhood of zero. Similarly one shows that
$v_h\to v_0$  in $C_{\rm loc}^1(\R^n)$.
Therefore, it is easily seen that
$v_0$ is a nontrivial positive solution of~\eqref{sollim}, that is $v_0\in G(z_0)$.
Since the map $\beta$ is strictly convex,
we can use Theorem~\ref{step1} by choosing in \eqref{prima}
$\varphi=0$ and
\begin{align*}
{\mathcal L}(x,s,\xi)&:=\alpha(z_0+\eps_hx)\beta(\xi)+
V(z_0+\eps_hx)\frac{s^p}{p} -K(z_0+\eps_h x)F(s),\\
h(x)=h_{\eps,k}(x)&:=(\underbrace{0,\dots,0}_{k-1},T(\eps x),
\underbrace{0,\dots,0}_{n-k}),\quad\text{for $\eps>0$ and $k=1,\dots,n$},
\end{align*}
for every $x\in\R^n$, $s\in\R^+$ and $\xi\in\R^n$,
the function $T\in C_c^1(\R^n)$ being chosen so that
$T(x)=1$ for $|x|\leq 1$, and $T(x)=0$ for $|x|\geq 2$.
In particular, $h_{\eps,k}\in C^1_c(\R^n,\R^n)$ and
\begin{equation*}
\partial_ih^j_{\eps,k}(x)=\eps \partial_iT(\eps x)\delta_{kj},\quad
\text{for every $x\in\R^n$, $\eps>0$, and $i,j$, and $k$}.
\end{equation*}
Then, it follows from~\eqref{prima} that
\begin{align*}
&0=\sum_{i=1}^n\int_{\R^n}
\eps\partial_iT(\eps x)\alpha(z_0+\eps_h x)
\partial_{\xi_i}\beta(\nabla v_h)\partial_kv_h\,dx \\
&-\int_{\R^n}\eps\partial_kT(\eps x)\Big[\alpha(z_0+\eps_h x)
\beta(\nabla v_h)+V(z_0+\eps_hx)
\frac{v_h^p}{p}-K(z_0+\eps_hx)F(v_h)\Big]\,dx\\
&-\int_{\R^n}\eps_hT(\eps x)\!\Big[\frac{\partial\alpha}{\partial
x_k}(z_0+\eps_h x)\beta(\nabla v_h)+\frac{\partial V}{\partial x_k}
(z_0+\eps_hx)\frac{v_h^p}{p}
-\frac{\partial K}{\partial x_k}
(z_0+\eps_hx)F(v_h)\Big]\,dx
\end{align*}
for every $\eps>0$, $h\in\N$, and $k=1,\dots,n$. Since the sequence $(v_h)$
is bounded in $W_V(\R^n)$,
by~\eqref{betacont},~\eqref{betacontxi}
and the boundedness of $\alpha$ and $K$, we have
\begin{gather*}
\left|\sum_{i=1}^n\int_{\R^n}
\partial_iT(\eps x)\alpha(z_0+\eps_h x)
\partial_{\xi_i}\beta(\nabla v_h)\partial_kv_h\,dx\right|\leq C, \\
\bigg|\int_{\R^n}\partial_kT(\eps x)\bigg[\alpha(z_0+\eps_h x)
\beta(\nabla v_h)+V(z_0+\eps_hx)
\frac{v_h^p}{p}-K(z_0+\eps_hx)F(v_h)\bigg]dx\bigg|\leq C',
\end{gather*}
for some positive constants $C,C'$. Therefore,
letting first $\eps\to 0$ yields
\begin{align}
\label{ps-prima}
\int_{\R^n}\bigg[\frac{\partial\alpha}{\partial x_k}(z_0+\eps_h x)
\beta(\nabla v_h)&+ \frac{\partial V}{\partial x_k}
(z_0+\eps_hx)\frac{v_h^p}{p}  \\
&-\frac{\partial K}{\partial x_k} \notag
(z_0+\eps_hx) F(v_h)\bigg]\,dx =0,
\end{align}
for every $h\in\N$, and $k=1,\dots,n$. Letting now $h \to \infty$,
by~\eqref{concentration}, we find
\begin{equation*}
\frac{\partial\alpha}{\partial x_k}(z_0)
\int_{\R^n}\beta(\nabla v_0)\,dx+\frac{\partial V}{\partial
x_k}(z_0)\int_{\R^n}\frac{v_0^p}{p}\,dx-\frac{\partial K}{\partial
x_k}(z_0)\int_{\R^n}F(v_0)\,dx=0,
\end{equation*}
for every $k=1,\dots,n$, which yields
\begin{equation*}
\nabla\alpha(z_0)\cdot w\int_{\R^n}\beta(\nabla v_0)\,dx
+\nabla V(z_0)\cdot w\int_{\R^n}\frac{v_0^p}{p}\,dx
=\nabla K(z_0)\cdot w\int_{\R^n}F(v_0)\,dx,
\end{equation*}
for every $w\in\R^n$.
Then, since $v_0\not\equiv 0$, assertion
{\rm (a)} immediately follows. Moreover, since $v_0\in G(z_0)$,
by the definition of $\partial\Gamma^-$, we obtain
\begin{align*}
\partial\Gamma^-(z_0;w)&=
\sup_{v\in G(z_0)}\bigg[\nabla\alpha(z_0)\cdot w\int_{\R^n}
\beta(\nabla v)\,dx\\
&+\nabla V(z_0)\cdot w
\int_{\R^n}\frac{|v|^p}{p}\,dx-\nabla K(z_0)\cdot w
\int_{\R^n}F(v)\,dx\bigg] \\
&\geq\nabla\alpha(z_0)\cdot w\int_{\R^n}
\beta(\nabla v_0)\,dx\\
&+\nabla V(z_0)\cdot w
\int_{\R^n}\frac{v_0^p}{p}\,dx-\nabla K(z_0)\cdot w
\int_{\R^n}F(v_0)\,dx=0,
\end{align*}
for every $w\in\R^n$.
Analogously, by the definition of $\partial\Gamma^+$,
we have
\begin{align*}
\partial\Gamma^+(z_0;w)&=
\inf_{v\in G(z_0)}\bigg[\nabla\alpha(z_0)\cdot w\int_{\R^n}
\beta(\nabla v)\,dx\\
&+\nabla V(z_0)\cdot w
\int_{\R^n}\frac{|v|^p}{p}\,dx-\nabla K(z_0)\cdot w
\int_{\R^n}F(v)\,dx\bigg] \\
&\leq\nabla\alpha(z_0)\cdot w\int_{\R^n}
\beta(\nabla v_0)\,dx\\
&+\nabla V(z_0)\cdot w
\int_{\R^n}\frac{v_0^p}{p}\,dx-\nabla K(z_0)\cdot w
\int_{\R^n}F(v_0) \,dx=0,
\end{align*}
for every $w\in\R^n$. Therefore
$z_0\in \mathfrak{C}$ and
assertion {\rm (b)} is proved. If $G(z_0)=\{v_0\}$, then clearly
$\Sigma$ admits all the directional derivatives at $z_0$, and
\begin{equation*}
\frac{\partial \Sigma}{\partial w}(z_0)=\partial\Gamma^-(z_0;w)=
\partial\Gamma^+(z_0;w)=0,\quad\text{for every $w\in\R^n$},
\end{equation*}
by virtue of {\rm (b)}. This proves assertion {\rm (c)}.
\end{proof}

The strong convergence required by \eqref{concentration} allows us to
take the limit as $h\to \infty$ in equation \eqref{ps-prima}. In the
semi-linear case one can construct uniform exponential
barriers for the family $(v_h)$, and therefore the strong
convergence of $(v_h)$ follows easily from the Lebesgue Convergence
Theorem (see~\cite{pompsec,wang,wangzeng}). The well-known loss of
regularity for solutions of quasi-linear equations is usually an
obstruction to this kind of argument. However, if the solutions belong
to a suitable space, then a pointwise concentration
suffices (see Corollary~\ref{stuart}).

\begin{remark}
We wish to point out that Theorem~\ref{mainth} holds true also
for the more general class of quasi-linear equations
\begin{equation*}
-\eps^p\dvg(\alpha(x)\partial_\xi\beta(u,\nabla u))
+\eps^p\alpha(x)\partial_s\beta(u,\nabla u)
+V(x)u^{p-1}=K(x)f(u),
\end{equation*}
under suitable assumptions on $\partial_\xi\beta(s,\xi)$
and $\partial_s\beta(s,\xi)$ (see~\cite{giacsqu}). On the other hand,
although the ground-state function $\Sigma$ can be
defined exactly as in Definition~\ref{groundef}
and $\Sigma(z)=c_z$ (cf.\ \cite[Theorem 3.2]{giacsqu}), the presence
of $u$ itself in the function $\beta$ makes the problems
of the regularity of $\Sigma$ and of the decay at infinity
for the rescaled family of solutions very complicated,
even in the nondegenerate case $p=2$.
\end{remark}

\begin{definition}\label{concentrationdef}
Let $z_0\in\R^n$. We say that a sequence $(u_{\eps_h})$
of solutions of~\eqref{problem} concentrates at $z_0$
if $u_{\eps_h}(z_0)\geq\ell>0$ for some $\ell>0$ and
for every $\eta>0$ there exist $\varrho>0$ and $h_0\in\N$ such that
\begin{equation*}
u_{\eps_h}(x)\leq \eta,\quad
\text{for every $h\geq h_0$ and $|x-z_0|\geq\eps_h\varrho$}.
\end{equation*}
\end{definition}
\noindent
This is precisely the notion of concentration
adopted in~\cite{wang,wangzeng}.

\begin{corollary}
\label{stuart}
Let $(u_{\eps_h})$ be a family of
solutions of~\eqref{problem} which concentrates
at a point $z_0\in\R^n$. Suppose that,
for every $h\in\N$ sufficiently large,
\[
u_{\eps_h} \in C_d^{1}(\R^n) \cap W^{2,n}(\R^n),
\]
where
\[
C_d^{1}(\R^n):=\Bigl\{ u\in C^1 (\R^n) :\,\,\,
\lim_{|x|\to \infty} u(x)=0\,\,\,
\text{and}\,\,\,\lim_{|x|\to \infty}
\nabla u(x)=0 \Bigr\}.
\]
Then, all the conclusions of~Theorem~\ref{mainth} hold true.
\end{corollary}
\begin{proof}
If $u_{\eps_h}\in C_d^{1}(\R^n)\cap W^{2,n}(\R^n)$, then one can
apply the results contained in~\cite{rs} to show
that the rescaled sequence $v_{\eps_h}$
decays exponentially fast at infinity,
uniformly with respect to $h$, together with all its partial
derivatives. Hence we can pass to the limit in
equation~\eqref{ps-prima}, and complete the
proof as in Theorem~\ref{mainth}.
\end{proof}

\noindent
For the particular, but important, case
$\alpha(x)=1$, $\beta(\xi)=|\xi|^{p-2}\xi$,
and $f(s)=s^{q-1}$, $p<q<p^*$, we can still prove
a fast-decay at infinity for the solutions.
\begin{lemma}
\label{decay}
Let $(u_{\eps_h})$ be a sequence of solutions of the problem
\begin{equation*}
\begin{cases}
-\eps^p \Delta_p u + V(x) u^{p-1} = K(x) u^{q-1} & \text{in $\R^n$}\\
u>0 & \text{in $\R^n$}
\end{cases}
\end{equation*}
which concentrates at $z_0\in\R^n$.
Then, if we set
$$
v_h(x):=u_{\eps_h}(z_0+\eps_hx),
$$
for each $\eta>0$ there
exist $R_\eta,C_\eta>0$ independent of $h$
such that
\begin{equation*}
|v_h(x)|\leq C_\eta\exp\bigg
\{-\left(\frac{\eta}{p-1}\right)^{1/p}|x|\bigg\},
\end{equation*}
for every $|x|\geq R_\eta$ and every $h\in\N$.
\end{lemma}
\begin{proof}
For every $h\in\N$, $v_h$ satisfies the equation
\begin{equation*}
-\Delta_pv_h+V(z_0+\eps_hx)v_h^{p-1}=K(z_0+\eps_hx)v_h^{q-1}\quad
\text{in $\R^n$}.
\end{equation*}
Since $(u_{\eps_h})$ is a concentrating sequence, it results that
\begin{equation*}
\lim_{|x|\to\infty}v_h(x)=0,\quad\text{uniformly in $h\in\N$}.
\end{equation*}
Then, setting $\inf_{x\in\R^n}V(x)=V_0$, given $\eta>0$ there
exists a positive constant $R_\eta$ independent of $h$ such that
\begin{equation*}
V(z_0+\eps_h x)v_h^{p-1}(x)-K(z_0+\eps_hx)v_h^{q-1}(x)
\geq (V_0-\eta)v_h^{p-1}(x),
\end{equation*}
for every $|x|\geq R_\eta$. It follows that the inequality
\begin{equation}
\label{first}
-\dvg(|\nabla v_h|^{p-2}\nabla
v_h)+(V_0-\eta)v_h^{p-1}\leq 0
\end{equation}
holds true for every $h\in\N$, and $|x|\geq R_\eta$.
Define now the function
\begin{equation*}
\Phi(x):=C_\eta\exp\left\{-\left(\frac{V_0-\eta}{p-1}\right)^{1
/p}|x|\right\},
\end{equation*}
where
\begin{equation*}
C_\eta:=\exp\left\{\left(\frac{V_0-
\eta}{p-1}\right)^{1/p}R_\eta\right\}\max_{|x|=R_\eta}v_h(x).
\end{equation*}
Notice that, since $v_h$ is uniformly bounded, we can assume
that $C_\eta$ is independent of $h$. Now, exactly the same
computations of~\cite[Theorem 2.8]{kabe} entail
\begin{equation}
\label{second}
-\dvg(|\nabla \Phi|^{p-2}\nabla \Phi)+(V_0-\eta)\Phi^{p-1}\geq 0.
\end{equation}
Testing inequalities~\eqref{first} and~\eqref{second}
with $\phi=(v_h-\Phi)^+$ on $\{|x|\geq R_\eta\}$ yields
\begin{align*}
&\int_{\{|x|\geq R_\eta\}\cap\{v_h>\Phi\}}\left(
|\nabla v_h|^{p-2}\nabla v_h\cdot\nabla (v_h-\Phi)+
(V_0-\eta)v_h^{p-1}(v_h-\Phi)\right)\,dx\leq 0, \\
&\int_{\{|x|\geq R_\eta\}\cap\{v_h>\Phi\}}\left(
|\nabla \Phi|^{p-2}\nabla \Phi\cdot\nabla (v_h-\Phi)+
(V_0-\eta)\Phi^{p-1}(v_h-\Phi)\right)\,dx\geq 0.
\end{align*}
By subtracting the previous inequalities,
and taking into account that
\begin{equation*}
\sum_{i=1}^n(|\xi|^{p-2}\xi_i-|\zeta|^{p-2}\zeta_i)
(\xi_i-\zeta_i)>0,\quad\text{for
every $\xi,\zeta\in\R^n$, $\xi\neq\zeta$},
\end{equation*}
we get
\begin{equation*}
\int_{\{|x|\geq R_\eta\}\cap\{v_h>\Phi\}}
(v_h^{p-1}-\Phi^{p-1})(v_h-\Phi)\,dx\leq 0.
\end{equation*}
Since $v_h$ and $\Phi$ are continuous functions, it has to be
\begin{equation*}
\{|x|\geq R_\eta\}\cap\{v_h>\Phi\}=\emptyset,
\quad\text{for every $h\in\N$},
\end{equation*}
which implies the assertion.
\end{proof}

\begin{theorem}
\label{plaplaci}
Let $(u_{\eps_h})$ be a
sequence of solutions of the problem
\begin{equation}
\label{plaplac}
\begin{cases}
-\eps^p\Delta_pu+V(x)u^{p-1}=K(x)u^{q-1} & \text{in $\R^n$}\\
u>0 & \text{in $\R^n$}
\end{cases}
\end{equation}
which concentrates at $z_0\in\R^n$. Then, the following facts hold:
\begin{itemize}
\item[(a)] the vectors $\nabla V(z_0)$ and
$\nabla K(z_0)$ are proportional\,;
\vskip3pt
\item[(b)] $z_0\in \mathfrak{C}$, that is
$z_0$ is a weak-concentration point for~\eqref{plaplac};
\vskip3pt
\item[(c)] if $1<p\leq 2$ then
all the partial derivatives of $\Sigma$ at $z_0$ exist, and
$\nabla\Sigma(z_0)=0$, that is $z_0$ is a critical
point of $\Sigma$.
\end{itemize}
\end{theorem}
\begin{proof}
By virtue of Lemma~\ref{decay}
we can pass to the limit in equation \eqref{ps-prima}
and get assertions {\rm (a)} and {\rm (b)} as in Theorem~\ref{mainth}.
If $1<p\leq 2$, by combining the results of~\cite{dr,gbli}
and~\cite{sertang}, for every $z\in\R^n$, problem~\eqref{limprob}
admits a unique positive $C^1$ solution
(up to translations) such that $u (x) \to 0$ as $|x| \to \infty$.
Then $G(z_0)=\{v_0\}$, and assertion {\rm (c)} follows by the
corresponding assertion in Theorem~\ref{mainth}.
\end{proof}

\section{The Semi-linear Case}

The main goal of this section is that of getting, in the particular case
$\beta(\xi)=\xi$, namely semi-linear equations, a more accurate
version of Theorem~\ref{mainth} involving the Clarke subdifferential
of the ground-state function $\Sigma$. We wish to stress that we have in mind the
case when $f$ is {\em not} simply the power nonlinearity $u^{p-1}$
(cf.\ Remark~\ref{smooth-power}).
\par
For $z\in\R^n$ fixed, we consider the limiting
functional $I_z:H^1(\R^n)\to\R$,
\begin{equation*}
I_{z}(u):=\alpha(z)\int_{\R^n}|\nabla u|^2\,dx
+\frac{V(z)}{p}\int_{\R^n}|u|^p\,dx-K(z)\int_{\R^n}F(u)\,dx
\end{equation*}
whose critical points are of course solutions of~\eqref{limprob}.
The minimax levels $c_z$ of $I_z$ are defined according to~\eqref{MPlev}.
Throughout the rest of this section, we will denote by
$S(z)$ the set of all the nontrivial solutions of~$(P_z)$ corresponding
to the energy level $\Sigma(z)$ (the set of ground-states).
It is known that $S(z)\neq\emptyset$ for every $z\in\R^n$ (see~\cite{berlio}).

\vskip2pt
\noindent
As the next lemma shows, in this particular
situation, the function $\Sigma$
has further regularity properties (and in some cases it relates to
the maps $\partial\Gamma^-$ and $\partial\Gamma^+$).

\begin{lemma}
\label{pro-sigma0}
If $p=2$ and condition~\eqref{increasing} holds,
then the following facts hold:
\begin{itemize}
\item[(i)] $\Sigma$ is locally Lipschitz;
\vskip3pt
\item[(ii)] the directional derivatives from the left and the
right of $\Sigma$ at $z$ along $w$, $\big(\frac{\partial\Sigma}{\partial w}\big)^{-}\!(z)$
and $\big(\frac{\partial\Sigma}{\partial w}\big)^{+}\!(z)$ respectively,
exist at every point $z\in\R^n$, and it holds
\begin{align*}
\left(\frac{\partial\Sigma}{\partial w}\right)^{-}\!\!(z)&=
\sup_{v\in S(z)}\nabla_z I_z(v)\cdot w, \\
\left(\frac{\partial\Sigma}{\partial w}\right)^{+}\!\!(z)&=
\inf_{v\in S(z)}\nabla_z I_z(v)\cdot w,
\end{align*}
for every $z,w\in\R^n$. In particular, if $G(z)=S(z)$, we have
\begin{equation}
\label{yformula}
\partial\Gamma^-(z;w)=\left(\frac{\partial
\Sigma}{\partial w}\right)^{-}\!\!(z)
\quad\text{and}\quad
\partial\Gamma^+(z;w)=\left(\frac{\partial
\Sigma}{\partial w}\right)^{+}\!\!(z),
\end{equation}
for every $w\in\R^n$.
\end{itemize}
\end{lemma}
\begin{proof}
By the results of~\cite{wangzeng},
$\Sigma$ is a locally Lipschitz map. We remark here
that, since $z$ acts as a parameter, the functional
$I_z$ is invariant under orthogonal change of variables.
Therefore, without loss of generality, to get the
formulas for the left and right directional
derivatives of $\Sigma$, it suffices to show that
\begin{align*}
& \left(\frac{\partial\Sigma}{\partial z_i}\right)^{-}\!\!(z)=
\sup_{v\in S(z)}\bigg[\frac{\partial\alpha}{\partial z_i}(z)\int_{\R^n}
\frac{|\nabla v|^2}{2}
+\frac{\partial V}{\partial z_i}(z)\int_{\R^n}\frac{|v|^p}{p}-
\frac{\partial K}{\partial z_i}(z)\int_{\R^n}F(v)\bigg],\\
& \left(\frac{\partial\Sigma}{\partial z_i}\right)^{+}\!\!(z)=
\inf_{v\in S(z)}\bigg[\frac{\partial\alpha}{\partial z_i}(z)\int_{\R^n}
\frac{|\nabla v|^2}{2}
+\frac{\partial V}{\partial z_i}(z)\int_{\R^n}\frac{|v|^p}{p}-
\frac{\partial K}{\partial z_i}(z)\int_{\R^n}F(v)\bigg],
\end{align*}
for every $z\in\R^n$ and $i=1,\dots,n$. These can be obtained
arguing as in~\cite{pompsec,wangzeng}. Finally, formulas~\eqref{yformula}
follow by the definition of $\partial\Gamma^+(z;w)$ and
$\partial\Gamma^-(z;w)$.
\end{proof}

\begin{remark}
\label{smooth-power}
Assume that $p=2$, $K$ is bounded from below
away from zero, and $f(u)=u^{q-1}$, where $2<q<2^*$. Then $\Sigma$ is smooth and it can
be given an explicit form (cf.~\cite[Remark 2.1]{pompsec}): there exists $C_q>0$ such that
\begin{equation*}
\Sigma(z)=C_q\left[\frac{V(z)}{K(z)}\right]^{\frac{q}{q-2}-\frac{n}{2}}
\!\!\!\!\sqrt{\alpha(z)K(z)},\quad\text{for every $z\in\R^n$}.
\end{equation*}
\end{remark}

Let us now recall from~\cite{clarke} two
definitions that will be useful in the sequel.
\begin{definition}
Let $f:\R^n\to\R$ be a locally Lipschitz function
near a given point $z\in\R^n$. The generalized derivative of the
function $f$ at $z$ along the direction $w\in\R^n$ is defined by
\begin{equation*}
f^0(z;w):=\limsup_{\substack{\xi \to z
\\ \lambda \to 0+}}\frac{f(\xi+\lambda w)-f(\xi)}{\lambda}.
\end{equation*}
\end{definition}

\begin{definition}
Let $f:\R^n\to\R$ be a locally Lipschitz function
near a given point $z\in\R^n$. The Clarke
subdifferential (or generalized gradient) of $f$ at $z$
is defined by
\begin{equation*}
\partial f(z):=\Big\{\eta\in\R^n:\,\,
f^0(z,w)\geq \eta\cdot w,\,\,\,\text{for every $w\in\R^n$}\Big\}.
\end{equation*}
\end{definition}

By~\cite[Proposition 2.3.1]{clarke} we learn that
\begin{proposition}
\label{oppositesub}
For every $z\in\R^n$, the set $\partial f(z)$ is
nonempty and convex, and
\begin{equation*}
\partial (-f)(z)=-\partial f(z).
\end{equation*}
\end{proposition}

The next is the main result of this section.
\begin{theorem}
\label{laplaci}
Assume that $(u_{\eps_h})$ is a
sequence of solutions of the problem
\begin{equation}
\label{problem0}
\begin{cases}
-\eps^2\dvg(\alpha(x)\nabla u)
+V(x)u=K(x)f(u) & \text{in $\R^n$}\\
u>0 & \text{in $\R^n$}
\end{cases}
\end{equation}
which concentrates at $z_0$. Then, the following facts hold:
\begin{itemize}
\item[(a)] the vectors
\begin{equation*}
\nabla\alpha(z_0),\,\,\,
\nabla V(z_0),\,\,\,
\nabla K(z_0)
\end{equation*}
are linearly dependent\,;
\vskip3pt
\item[(b)] $z_0\in\mathfrak{C}$, that is
$z_0$ is a weak-concentration point for~\eqref{problem0};
\vskip3pt
\item[(c)] if either $G(z_0)=S(z_0)$ or
\begin{equation}
\label{limc}
\eps_h^{-n} J_{\eps_h}(u_{\eps_h}) \to c_{z_0},
\end{equation}
where
\begin{align}
\label{funzio}
J_{\eps}(v)&=\frac{\eps^2}{2}
\int_{\R^n} \alpha (x)|\nabla v|^2 \, dx + \frac{1}{2}
\int_{\R^n} V(x) |v|^2 \, dx
- \int_{\R^n} K(x) F(v)\, dx,
\end{align}
we have
\begin{equation*}
0\in\partial\Sigma(z_0),
\end{equation*}
that is $z_0$ is a critical point of $\Sigma$
in the sense of the Clarke subdifferential\,;
\vskip3pt
\item[(d)] if $S(z_0)=\{v_0\}$, then
all the partial derivatives of $\Sigma$ at $z_0$ exist, and
\begin{equation*}
\nabla\Sigma(z_0)=0,
\end{equation*}
that is $z_0$ is a critical point of $\Sigma$.
\end{itemize}
\end{theorem}
\begin{proof}
For problem~\eqref{problem0} it is possible to
prove the existence of uniform exponentially decaying
barriers. Then we can pass to the limit in equation \eqref{ps-prima},
to get assertions {\rm (a)} and {\rm (b)} as
in Theorem~\ref{mainth}.
If $G(z_0)=S(z_0)$, by combining
formulas~\eqref{yformula} of Lemma~\ref{pro-sigma0}
with {\rm (b)} of Theorem~\ref{mainth}, we have
\begin{equation}
\label{frdirectder}
\left(\frac{\partial\Sigma}{\partial w}\right)^-(z_0)\geq 0
\quad\text{and}\quad
\left(\frac{\partial\Sigma}{\partial w}\right)^+(z_0)\leq 0,
\end{equation}
for every $w\in\R^n$. In particular, it holds
\begin{equation*}
\left(\frac{\partial (-\Sigma)}{\partial w}\right)^{+}\!\!(z_0)\geq
0,\quad\text{for every $w\in\R^n$}.
\end{equation*}
Then, by the definition of $(-\Sigma)^0(z_0;w)$ we get
\begin{equation*}
(-\Sigma)^0(z_0;w)\geq\left(\frac{\partial (-\Sigma)}
{\partial w}\right)^{+}\!\!(z_0)\geq 0,\quad
\text{for every $w\in\R^n$}.
\end{equation*}
By the definition of $\partial(-\Sigma)(z_0)$ we immediately get
$0\in \partial(-\Sigma)(z_0)$, which, together with Proposition~\ref{oppositesub},
yields assertion {\rm (c)}. To prove the same conclusion
when~\eqref{limc} holds, we simply remark that $c_{z_0}=\Sigma (z_0)$.
Therefore, if $v_0$ is the limit of
the sequence $(v_h)$ defined in \eqref{vh},
then $v_0 \in S(z_0)$ because we
can exploit again some exponential
barrier to pass to the limit.
As a consequence, arguing as in Theorem~\ref{mainth}, it follows that
inequalities~\eqref{frdirectder} hold
and we are reduced to the previous case.
Finally, if $S(z_0)=\{v_0\}$,
the map $\Sigma$ admits all the directional
derivatives at $z_0$ and, by virtue
of~\eqref{frdirectder} they are equal to zero,
which proves (d).
\end{proof}

We would like to remark that a different
definition of concentration has been
used in~\cite{gm}. We recall it here,
suitably adapted to our purposes.

\begin{definition}
Assume that $u_\eps \in C^2 (\R^n)$ is a family of
solutions of~\eqref{problem0} and let $J_\eps$ be as in~\eqref{funzio}.
Moreover, let $x_\eps\in\R^n$ be such that $\max_{x\in\R^n} u_\eps = u_\eps (x_\eps).$
We say that $u_\eps$ concentrates at $z_0\in\R^n$ if
the following facts hold:
\begin{itemize}
\item[(i)] $\lim\limits_{\eps\to 0}\,x_\eps=z_0$\,;
\vskip2pt
\item[(ii)] $\lim\limits_{\eps\to 0}\,\eps^{-n}J_\eps(u_\eps)=c_{z_0}$.
\end{itemize}
\end{definition}

It is not difficult to check that if $(u_\eps)$ is a sequence as in the
above definition, then $(u_\eps)$ concentrates at $z_0$ in the sense of
Definition~\ref{concentrationdef}, vanishing at an exponential
rate away from $z_0$ (cf.\ \cite[Lemma 4.2]{gm}). In particular,
according to {\rm (c)} of Theorem~\ref{laplaci}, we have
$0\in\partial\Sigma(z_0)$.

\vskip8pt
\noindent
We finish the paper with an open problem. Assume that $(u_h)$
is a sequence of solutions of problem~\eqref{problem0}. Suppose that
these solutions concentrate at $z_0 \in \R^n$, and $S(z_0)=\{v_0\}$.
Is it possible to prove that $z_0$ is a $C^1$-stable critical
point of $\Sigma$, according to the
definition of Yanyan Li~\cite{yyli}?

\medskip
\subsection*{Acknowledgment}
The authors are indebted to the anonymous referee for her/his
careful reading of the manuscript and for valuable remarks and comments.

\end{document}